\documentclass[12pt]{article}
\usepackage{latexsym,amsmath,graphicx}

\newtheorem{lemma}{Lemma}[section]
\newtheorem{theorem}[lemma]{Theorem}
\newtheorem{proposition}[lemma]{Proposition}

\newcommand {\qed}{\rule{2mm}{2mm}}

\newtheorem{example}[lemma]{Example}{\normalfont}
\let\olddefinition\example
\renewcommand{\example}{\olddefinition\normalfont}

\begin{document}
\title{On the upper embedding of Steiner triple systems and Latin squares}
\author{Terry S. Griggs\\
School of Mathematics and Statistics\\
The Open University\\
Milton Keynes MK7 6AA\\
UNITED KINGDOM\\
{\tt t.s.griggs@open.ac.uk}\\~\\
Thomas A. McCourt\\
Craigslea State High School\\
Brisbane QLD 4032\\
AUSTRALIA\\
{\tt tom.a.mccourt@gmail.com}\\~\\
Jozef \v Sir\'a\v n\\
School of Mathematics and Statistics\\
The Open University\\
Milton Keynes MK7 6AA\\
UNITED KINGDOM\\
{\tt j.siran@open.ac.uk}\\~\\
and Slovak University of Technology\\
Bratislava 81005\\
SLOVAKIA\\~\\}
\date{} \maketitle

\newpage
\begin{abstract}
\noindent It is proved that for any prescribed orientation of the triples of either a Steiner triple system or a Latin square of odd order,
there exists an embedding in an orientable surface with the triples forming triangular
faces and one extra large face.
\end{abstract}

\bigskip


\newpage
\section{Introduction}\label{intro}
The motivation for the work described herein comes from two previous papers, respectively on the upper embedding of Steiner triple systems \cite{GGS} and the upper embedding of Latin squares \cite{GPS}. First we recall the relevant definitions. Let ${\cal X}=(V,{\cal B})$ be a (partial) triple system on a point set $V$, that is, a collection ${\cal B}$ of $3$-element subsets of $V$, called blocks or triples, such that every $2$-element subset of $V$ is contained in at most one triple in ${\cal B}$. Equivalently, such a triple system ${\cal X}$ may be viewed as a pair $(K,{\cal B})$, where $K$ is a graph with vertex set $V$ and edge set $E$ consisting of all pairs $uv$ for distinct points $u,v\in V$ such that $\{u,v\}$ is a subset of some block in ${\cal B}$. In other words, in the graph setting, ${\cal B}$ is regarded as a decomposition of the edge set $E$ of $K$ into triangles. Of course, such a graph $K=(V,E)$ may admit many decompositions into triangles and so the set of blocks ${\cal B}$ needs to be specified. We will refer to $K=(V,E)$ with the specified decomposition ${\cal B}$ of $E$ as the graph {\em associated} with the triple system ${\cal X}=(V,{\cal B})$. We will be assuming throughout that the triple systems considered here are {\em connected}, meaning that their associated graphs are connected. In the case
of a Steiner triple system ${\cal S}$ = STS($n$) where $n=|V|$, the associated graph is the complete graph $K_n$. Such systems exist if and only if $n \equiv$ 1 or 3 (mod 6) \cite{K}. For a Latin square L = LS($n$), the associated graph is the complete tripartite graph $K_{n,n,n}$ where the three parts of the tripartition are the rows, the columns and the entries of the Latin square.

By an embedding of a triple system ${\cal X}=(V,{\cal B})$ we will understand a cellular embedding $\vartheta:\ K\to \Sigma$ of the associated graph $K=(V,E)$ of $\cal X$ in an orientable surface $\Sigma$, such that every triangle in ${\cal B}$ bounds a face of $\vartheta$. Such faces will be called {\em block faces}, and the remaining faces of the embedding $\vartheta$ will be called {\em outer faces}. By the properties of the set ${\cal B}$, in the embedding $\vartheta$ every edge of $E$ lies on the boundary of exactly one block face, so that there is at least one outer face in $\vartheta$. The extreme case occurs if such an embedding has exactly one outer face; we then speak about an {\em upper embedding} and call the triple system ${\cal X}=(V,{\cal B})$ {\em upper embeddable}.

A necessary and sufficient condition for upper embeddability of triple systems follows from available knowledge about upper embeddings of graphs in general. To make use of this we will represent triple systems by their point-block incidence graphs as usual in design theory. For a triple system ${\cal X}=(V,{\cal B})$ its {\em point-block incidence graph} is the bipartite graph $G({\cal X})$ with vertex set $V\cup {\cal B}$ and edge set consisting of pairs $\{v,B\}$ for $v\in V$ and $B\in {\cal B}$ such that $v \in B$. The pair $(V,{\cal B})$ forms the bi-partition of the vertex set of $G({\cal X})$; vertices in $V$ and ${\cal B}$ will be referred to as {\em point vertices} and {\em block vertices} respectively. By our convention regarding triple systems, the graph $G({\cal X})$ is assumed to be connected, and note that every block vertex has valency $3$ in $G({\cal X})$. 

It is a folklore fact in topological design theory that there is a one-to-one correspondence between orientable embeddings of a triple system ${\cal X}$ and its point-block incidence graph $G({\cal X})$ in orientable surfaces; the correspondence is illustrated in Fig. 1. As is obvious from this figure, an embedding of $G({\cal X})$ arises from an embedding of ${\cal X}$ naturally. On the other hand, given an embedding of $G({\cal X})$ one obtains the corresponding embedding of ${\cal X}$ by `inflating' every block vertex into a triangle on the surface. In particular, a triple system ${\cal X}$ is upper embeddable if and only if its point-block incidence graph is embeddable with exactly one face; such graphs are also called upper-embeddable. By a classical result of Jungerman \cite{J} and Xuong \cite{X}, a graph (in particular, the point-block incidence graph of a triple system) is upper-embeddable if and only if the graph contains a spanning tree such that each of its co-tree components has an even number of edges.

\begin{figure}[h]
\begin{center}
      \scalebox{0.5}
       {\includegraphics{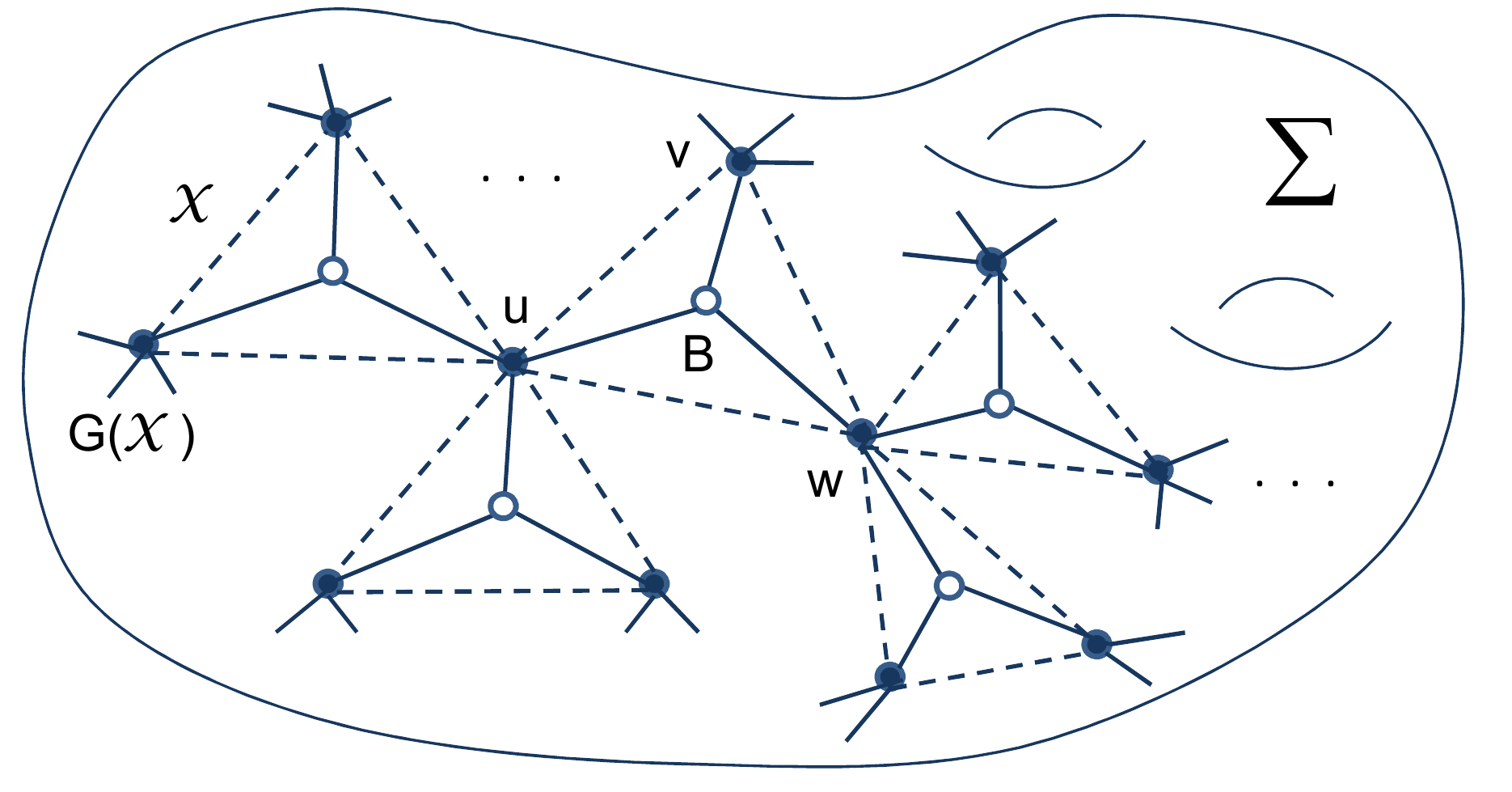}}
       \vskip -2mm
       \caption{Correspondence between embeddings of ${\cal X}$ and $G({\cal X})$.}
\end{center}
    \end{figure}

\vskip -5mm

Let us now look more closely at an upper embedding of a triple system ${\cal X}=(V,{\cal B})$ in an orientable surface $\Sigma$, or, equivalently, at an upper embedding $\vartheta:\ G({\cal X})\to \Sigma$ of the point-block intersection graph of ${\cal X}$ in $\Sigma$. For any triple $B=\{u,v,w\}\in {\cal B}$ the (given) orientation of $\Sigma$ induces one of the two cyclic permutations $(u,v,w)$ or $(u,w,v)$ when looking at the three points from the `centre' of the triangular face representing $B$. Equivalently for any block vertex $B$ of $G({\cal T})$ with point vertex neighbours $u$, $v$, $w$, the orientation of $\Sigma$ induces one of the cyclic permutations $(u,v,w)$ or $(u,w,v)$ of the three points `around' the block vertex $B$. We refer to any of the two permutations as an {\em orientation of the triple}, or equivalently, an {\em orientation of the neighbourhood}. In this terminology, the orientation of $\Sigma$ induces one of the two possible orientations of every triple (or, of the neighbourhood of every block vertex) in the embedding. It was proved in \cite{GGS} that every Steiner triple system STS($n$) and in \cite{GPS} that every Latin square LS($n$) where $n$ is odd has an upper embedding in an orientable surface. In the Latin square case the restriction that $n$ must be odd is determined by Euler's formula $(V+F-E=2-2g)$. However in both \cite{GGS} and \cite{GPS} no attention was paid to the orientation of the triples in the upper embeddings.

The driving question of our research is the following question. Which triple systems ${\cal X}=(V,{\cal B})$ have the property that, for every choice of the orientation of all triples $B\in {\cal B}$, the system ${\cal X}$ admits an upper embedding in an orientable surface such that the orientation of the surface induces the preassigned orientation of $B$ for {\em every} triple $B\in {\cal B}$? We will refer to this property simply as {\em upper embeddability in every orientation of triples.} The main results of this paper are the two theorems below which substantially extend the results in the two papers \cite{GGS} and \cite{GPS}.

\begin{theorem}\label{thm:STS}
Every Steiner triple system admits an upper embedding in every orientation of triples.
\end{theorem} 
\begin{theorem}\label{thm:LSodd}
Every Latin square of odd order admits an upper embedding in every orientation of triples.
\end{theorem}
 
\section{Spanning tree}\label{span}
In what follows we prove a sufficient condition for a triple system to admit upper embeddability in every orientation of triples. The result will be stated in terms of the point-block incidence graph of a triple system and, as one expects by Jungerman and Xuong's Theorem \cite{J} and \cite{X}, the statement will involve existence of a spanning tree with particular properties. We recall the connectivity assumption of our triple systems.

\begin{theorem}\label{thm:upper} 
Let ${\cal X}$ be a triple system and let $G=G({\cal X})$ be its point-block incidence graph. If $G$ admits a spanning tree such that every point vertex has even valency in the corresponding co-tree, then ${\cal X}$ admits an upper embedding in every orientation of triples.
\end{theorem}

\textbf{Proof.} Let ${\cal X}=(V,{\cal B})$ and suppose that every triple $B=\{u,v,w\}\in {\cal B}$ has been assigned an orientation, i.e., one of the two cyclic permutations $(u,v,w)$, $(u,w,v)$. Using the ideas of Jungerman and Xuong we will show (by induction) how to build an upper embedding of the point-block incidence graph $G=G({\cal X})$ of ${\cal X}$ in an oriented surface in such a way that its orientation will induce the preassigned orientation on every $B\in {\cal B}$.

Let $T$ be a spanning tree of $G$ as in the assumption of our theorem, that is, such that every point vertex of the subgraph of $G$ induced by the set $E(G){\setminus}E(T)$ of co-tree edges has even valency. Since $G$ is bipartite, this assumption implies that the set $E(G){\setminus}E(T)$ has a decomposition into paths of length two that have a point vertex in the centre; in particular, the number of co-tree edges here is even. We note that the quantity $\beta=|E(G){\setminus}E(T)|$ is known as the Betti number of $G$. The set of co-tree edges thus decomposes into $\beta/2$ paths $P$ of the form $P=AuB$, where $A,\, B$ are block vertices and $u$ is a point vertex.

To proceed, we prove quite a general auxiliary statement on upper embeddability of extensions of spanning subgraphs of our point-block incidence graph by paths as above. Let $H$ be a connected spanning subgraph of $G$ and let $A$ and $B$ be block vertices and $u$ be a point vertex of $H$ such that $P=AuB$ is a path in $G$ but the edges $Au$ and $Bu$ are not in $H$. Further, let ${\cal C}$ be the set of block  vertices of $G$ that have valency $3$ in $H$. We will say that $H$ {\em upper embeds in every orientation at ${\cal C}$} if, for every vertex $C\in {\cal C}$ and every orientation of the neighbourhood of $C$, there is an embedding of $H$ in an orientable surface such that its orientation induces the preassigned orientation around every vertex $C\in {\cal C}$. Extending this terminology in a natural way to the graph $H'=H\cup P$ and the set ${\cal C}'$ of block vertices of valency $3$ in $H'$, we prove the following.
\smallskip

\noindent{\bf Claim.} {\sl If $H$ upper embeds in every orientation at ${\cal C}$, then $H'$  upper embeds in every orientation at ${\cal C}'$.}
\smallskip

To prove our Claim, let $H\to \Sigma$ be an upper embedding of $H$ in an orientable surface $\Sigma$ such that its orientation induces the preassigned orientations of neighbourhoods of block vertices in ${\cal C}$. Let $P=AuB$ be a path as above. The boundary of the single face $F$ of this embedding is, without loss of generality, a closed walk in $H$ of the form $(uXAYBZ)$, where $X,\, Y,\, Z$ are $u\to A$, $A\to B$ and $B\to u$ walks of $H$ traversed in the direction induced by the (clockwise) orientation of $\Sigma$, as one can see in Fig. 2 when disregarding the arcs (edges with direction) $a,b,c,d$. We point out that in our considerations the order of appearance of the vertices $A$ and $B$ on the boundary of $F$ will be immaterial.

Ignoring the condition on the set ${\cal C}'$, the embedding of $H$ can be extended to an upper embedding of $H'$. Namely, letting $a=uA$ and $b=uB$ denote the arcs from $u$ to $A$ and $u$ to $B$ respectively, one `adds' the path $P=AuB$ to the single face $F$ of $\vartheta$ in such a way that all the local cyclic orderings of arcs emanating from vertices distinct from $A,\, u,\, B$ are kept intact and the local cyclic ordering of neighbours of $u$ is extended from $(\ldots, c,d,\ldots)$ to $(\ldots,c,b,a,d, \ldots)$ as in Fig. 2; the local cyclic order of arcs at $A$ and $B$ is obvious from Fig. 2. 

\begin{figure}[h]
\begin{center}
      \scalebox{0.55}
       {\includegraphics{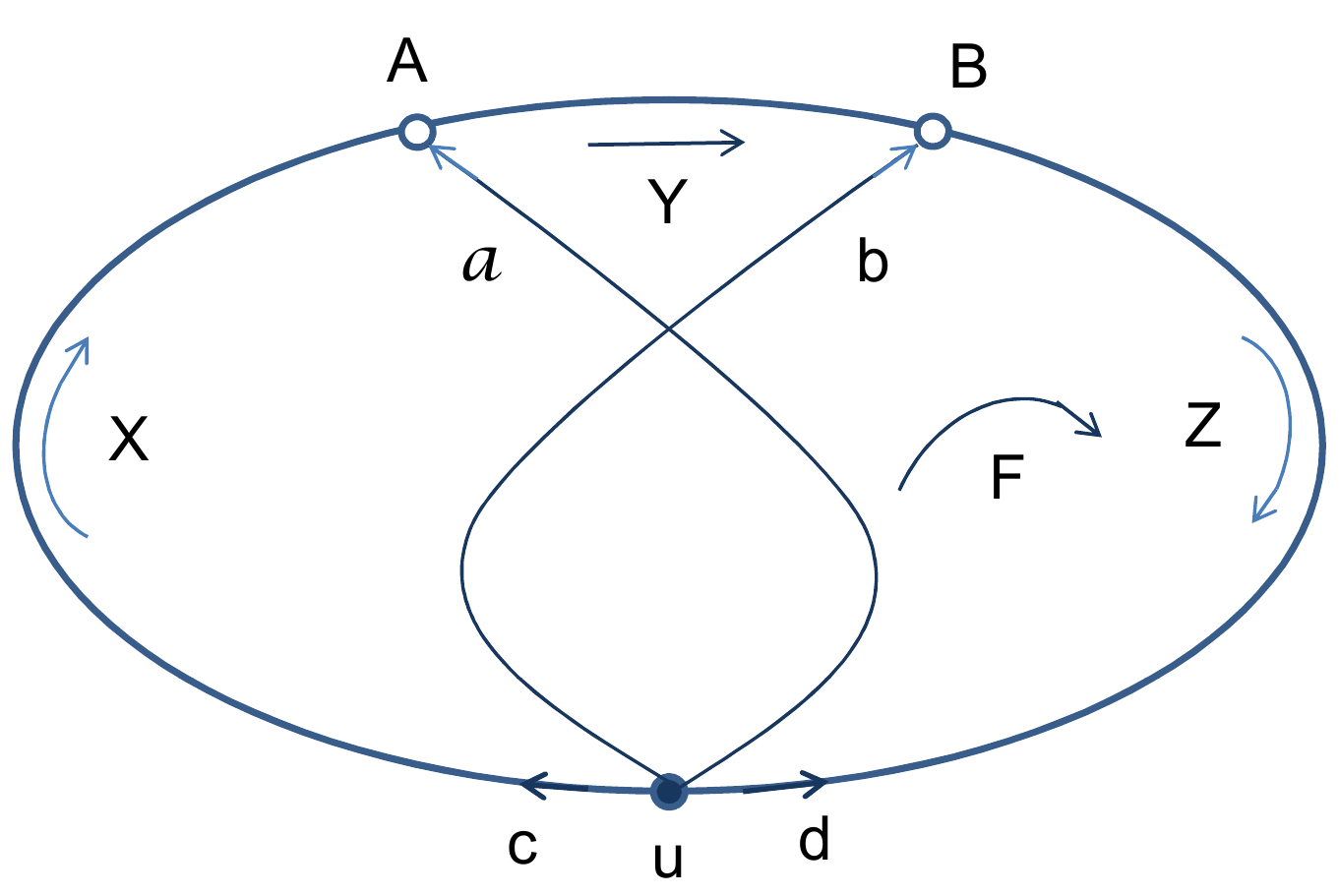}}
       \vskip -2mm
       \caption{Extending $F$ by adding the path $P=AuB$.}
\end{center}
    \end{figure}

\noindent The single face $F'$ of the new embedding of $H'$ in an oriented surface $\Sigma'$ is obtained by tracing down its boundary, which is the closed walk $(uXa^-bZaYb^-)$, where $a^-$ and $b^-$ indicate traversals along $a$ and $b$ in the opposite direction. We again emphasize that, for this construction, the position of $A$ and $B$ on the boundary of $F$ is irrelevant. Note that the genus of $\Sigma'$ is equal to the genus of $\Sigma$ increased by $1$.

We now show that the above construction can be carried out in such a way that the orientation of $\Sigma'$ induces the preassigned orientations of neighbours of vertices in ${\cal C}'$. This orientation is obviously maintained for block vertices in the subset ${\cal C}\subset {\cal C}'$ and so one only needs to consider the block vertices in ${\cal C}'{\setminus}{\cal C}\subset \{A,B\}$; note that $A,B\notin {\cal C}$. If neither $A$ nor $B$ are in ${\cal C}'$, we simply use the above construction. If one or both of $A,\, B$ are in ${\cal C}'{\setminus}{\cal C}$, then we proceed as follows. 

Say, without loss of generality, that $A\in {\cal C}'$, which means that the valency of $A$ in $H$ must have been equal to $2$. We know that $u\in A$; let $A=\{u,v,w\}$ with $\{v,w\}$ being be the point vertices constituting the neighbourhood of $A$ in $H$. Since the boundary of the single face $F$ in $H$ must contain every edge twice (and traversed in each direction once) and $A$ has valency $2$ in $H$, it follows that the boundary of $F$ must have the form $(u\ldots vAw\ldots wAv\ldots)$, as displayed in Fig. 3; the relative position of the $vAw$ and $wAv$ paths is without loss of generality. 

\begin{figure}[h]
\begin{center}
      \scalebox{0.55}
       {\includegraphics{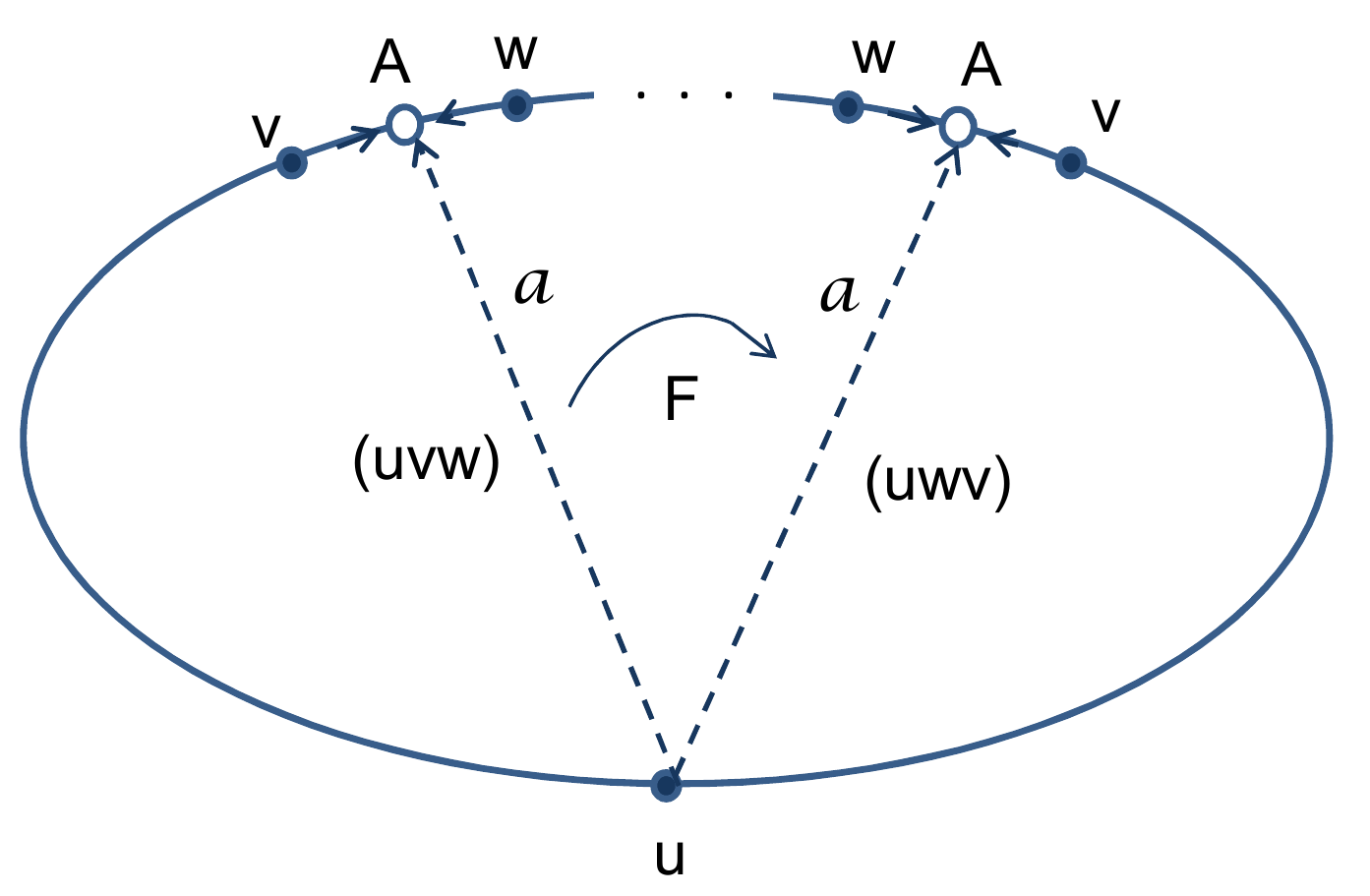}}
       \vskip -2mm
       \caption{Adding the edge $uA$ to control the orientation of the triple $A$.}
\end{center}
    \end{figure}

\vskip -5mm

Now, if the preassigned orientation of the block $A=\{u,v,w\}$ is given by the cyclic permutation $(uvw)$, we will add the arc $a=uA$ pointing to the position of $A$ appearing on the left-hand side of Fig. 3; if the preassigned orientation is $(uwv)$ we use the arc $a$ pointing to the occurrence of $A$ on the right-hand side in Fig. 3. If the vertex $B$ also has valency $2$ in $H$, we make a similar choice for the position of $B$ for addition of the arc $b=uB$. We then complete the construction of the single face embedding of $H'$ with the required properties as in the previous paragraph.

Having proved our Claim, the rest of the proof is straightforward. We begin by embedding the spanning tree $T$ in a sphere in such a way that its orientation induces a preassigned orientation at every block vertex that has valency $3$ in $T$. We then apply the construction of our Claim $\beta/2$ times for every path in the decomposition of the set $E(G){\setminus}E(T)$ of co-tree edges into $\beta/2$ paths whose middle vertex is a point vertex. As a result we obtain the upper embeddability of $G$ in every orientation at its set ${\cal B}$ of block vertices.
\hfill $\Box$

\section{Steiner triple systems}\label{STS}
We deal first with the case of Steiner triple systems. In view of the previous section, Theorem \ref{thm:STS} follows immediately from the following Proposition.
\begin{proposition}\label{spanSTS}
Let $S=(V,\mathcal{B})$ be a Steiner triple system of order $n$
and let $G=G(S)$ be its point-block incidence graph. Then $G$ admits a spanning tree such that every point vertex has even valency in the corresponding co-tree.
\end{proposition}
\textbf{Proof.}
Let $V=\{0,1,\ldots,n-1\}$. Construct a spanning tree $T$ of the
point-block incidence graph $G$ as follows. For convenience, we will refer to point vertices
simply as points and block vertices as blocks.  
Let the root (Level 0) of the tree be the point $0$. Connect the point to all $(n-1)/2$ blocks
containing it, which will be at Level 1. At Level 2 put all the  $n-1$ points $x \in V \setminus \{0\}$
and connect these to the blocks at Level 1 which contain them.
The remaining $n(n-1)/6-(n-1)/2=(n-1)(n-3)/6$ blocks are at Level 3.
Connect each one of these to a point at Level 2 which is contained in the block.
The structure of the spanning tree is shown below.

\begin{figure}[h]
\begin{center}
      \scalebox{0.45}
       {\includegraphics{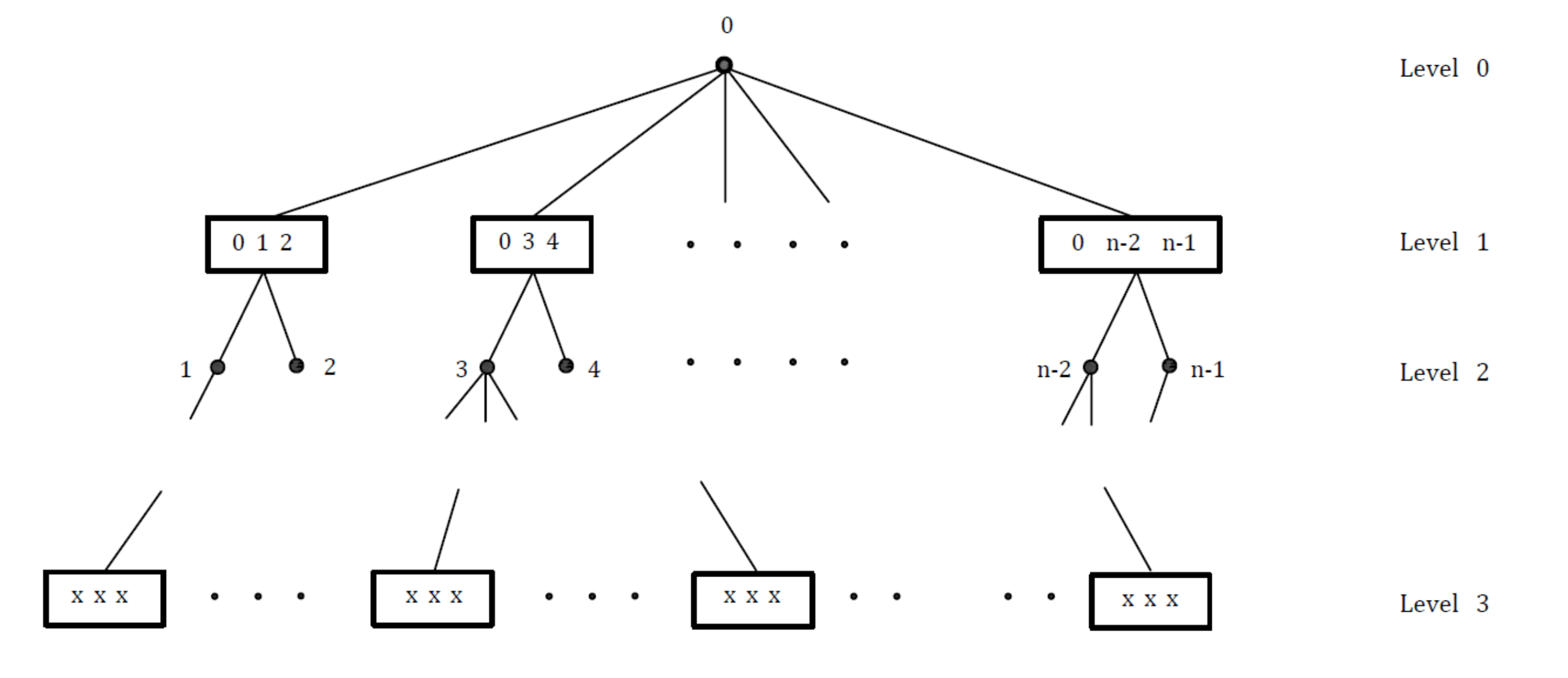}}
       \vskip -2mm
       \caption{Spanning tree of $G(S)$.}
\end{center}
    \end{figure}

\vskip -5mm

In the co-tree the point $0$ has zero valency and the points at level 2
have either odd or even valency including zero. We will refer to such points
as either \emph{odd points} or \emph{even points} respectively.
The number of edges in the co-tree is $3n(n-1)/6-n(n-1)/6-n+1=(n-1)(n-3)/3$
which is even. Thus the number of odd points is even.
The aim is to construct a spanning tree with no odd points.
The proof is in the form of an iterative algorithm.
Beginning with a spanning tree constructed as above, at any stage, if there are odd points
at Level 2, modify how the blocks at Level 3 are connected to the points at Level 2 in order to reduce the number of odd points until they are absent.
Let $P$ be the number of odd points. Proceed as follows.

Step 1. If $P=0$, then we are done; otherwise choose two odd points $a,b \in V$.

Step 2. Suppose that $\{0,a,b\} \in \mathcal{B}$; otherwise go to Step 3.
Then there exists a block $\{a,c,x\}$, $c,x,\in V \setminus \{0,a,b\}$ where
the edge $\{c,\{a,c,x\}\} \in T$ for otherwise the valency of the point $a$ in 
the co-tree would be zero. Consider the block $\{b,c,y\}$. There are three
possibilities to consider.\\ 
(i) The edge $\{c,\{b,c,y\}\} \in T$. Remove both $\{c,\{a,c,x\}\}$ and
$\{c,\{b,c,y\}\}$ from $T$ and replace by $\{a,\{a,c,x\}\}$ and $\{b,\{b,c,y\}\}$.
Return to Step 1.\\
(ii) The edge $\{b,\{b,c,y\}\} \in T$. Remove both $\{c,\{a,c,x\}\}$ and
$\{b,\{b,c,y\}\}$ from $T$ and replace by $\{a,\{a,c,x\}\}$ and $\{c,\{b,c,y\}\}$.
Return to Step 1.\\
(iii) The edge $\{y,\{b,c,y\}\} \in T$. Remove $\{c,\{a,c,x\}\}$
from $T$ and replace by $\{a,\{a,c,x\}\}$. If $c$ is an odd point, then both
$a$ and $c$ become even points and we return to Step 1. If $c$ is an even point,
then it becomes an odd point and $a$ becomes an even point. The block $\{b,c,y\}$
now contains the odd points $b$ and $c$ and $y \neq 0$ which may be odd or even.
If $y$ is an odd point then remove $\{y,\{b,c,y\}\}$ and replace by either
$\{b,\{b,c,y\}\}$ or $\{c,\{b,c,y\}\}$ and return to Step 1.
If $y$ is an even point then go to Step 4.

Step 3. We have a block $\{a,b,y\} \in \mathcal{B}$ where $y \neq 0$.
There are three possibilities to consider.\\
(i) The edge $\{a,\{a,b,y\}\} \in T$. Remove $\{a,\{a,b,y\}\}$ from $T$ and replace by
$\{b,\{a,b,y\}\}$. Return to Step 1.\\
(ii) The edge $\{b,\{a,b,y\}\} \in T$. Remove $\{b,\{a,b,y\}\}$ from $T$ and replace by
$\{a,\{a,b,y\}\}$. Return to Step 1.\\
(iii) The edge $\{y,\{a,b,y\}\} \in T$. If $y$ is an odd point, then remove $\{y,\{a,b,y\}\}$
from $T$  and replace by either $\{a,\{a,b,y\}\}$ or $\{b,\{a,b,y\}\}$. Return to Step 1.
If $y$ is an even point, then continue to Step 4.

Step 4. If we reach this step, we have a block $\{p,q,y\}$ where $p$ and $q$ are odd points,
$y$ is an even point and $\{y,\{p,q,y\}\} \in T$. (If we have come from Step 2,
$\{p,q\}=\{b,c\}$ and if we have come from Step 3, $\{p,q\}=\{a,b\}$.)

Now consider the odd point $p$. It is contained
in a block $\{0,p,p'\}$ to which it is connected by an edge
and the block $\{p,q,y\}$ to which it is not connected. There are $N=(n-5)/2$
further blocks containing $p$. Suppose that $f$ of these blocks $\{v_i,w_i,p\}$,
$i= 1,2,\ldots,f$ are not connected to $p$, i.e. $\{p,\{v_i,w_i,p\}\} \notin T$. So without loss
of generality let $\{v_i,\{v_i,w_i,p\}\} \in T$ if $i=1,2,\ldots,f$. Then the $N-f$ blocks
$\{v_i,w_i,p\}$, $i=f+1,f+2,\ldots,N$ are connected to $p$, i.e.
$\{p,\{v_i,w_i,p\}\} \in T$. Clearly $0 \leq f \leq N$ and since $p$ is an odd point,
then $f$ is even. Proceed as follows.\\
(A) If any of the points $v_i$, $i= 1,2,\ldots,f$ are odd points, say $v_\alpha$, then
remove the edge $\{v_\alpha,\{v_\alpha,w_\alpha,p\}\}$ and replace it by
$\{p,\{v_\alpha,w_\alpha,p\}\}$. Return to Step 1.\\ 
(B) If any of the points $v_i$ or $w_i$, $i= f+1,f+2,\ldots,N$ are odd points,
again say without loss of generality $v_\alpha$,
then remove the edge $\{p,\{v_\alpha,w_\alpha,p\}\}$
and replace it by $\{v_\alpha,\{v_\alpha,w_\alpha,p\}\}$. Return to Step 1.

If neither of the conditions (A) nor (B) is satisfied consider the odd point $q$
and repeat the above procedure. The point $q$ is contained
in a block $\{0,q,q'\}$ to which it is connected by an edge
and the block $\{p,q,y\}$ to which it is not connected. There are $N=(n-5)/2$
further blocks containing $q$. Suppose that $g$ of these blocks $\{v_i,w_i,q\}$,
$i= 1,2,\ldots,g$ are not connected to $q$, i.e. $\{q,\{v_i,w_i,q\}\} \notin T$. So without loss
of generality let $\{v_i,\{v_i,w_i,q\}\} \in T$ if $i=1,2,\ldots,g$. Then the $N-g$ blocks
$\{v_i,w_i,q\}$, $i=g+1,g+2,\ldots,N$ are connected to $q$, i.e.
$\{q,\{v_i,w_i,q\}\} \in T$. Clearly $0 \leq g \leq N$ and since $q$ is an odd point,
then $g$ is even. Proceed as follows.\\
(C) If any of the points $v_i$, $i= 1,2,\ldots,g$ are odd points, say $v_\beta$, then
remove the edge $\{v_\beta,\{v_\beta,w_\beta,q\}\}$
and replace it by $\{q,\{v_\beta,w_\beta,q\}\}$. Return to Step 1.\\ 
(D) If any of the points $v_i$ or $w_i$, $i= g+1,g+2,\ldots,N$ are odd points,
again say without loss of generality $v_\beta$,
then remove the edge $\{q,\{v_\beta,w_\beta,q\}\}$
and replace it by $\{v_\beta,\{v_\beta,w_\beta,q\}\}$.
Return to Step 1.
 
If neither of the conditions (C) nor (D) is satisfied consider the even point $y$.
Recall that $\{y,\{p,q,y\}\} \in T$. Remove this edge and replace it by either
$\{p,\{p,q,y\}\}$ or $\{q,\{p,q,y\}\}$. So now either $p$ or $q$ is an even point
and $y$ is an odd point. Repeat the above procedure for the point $y$. It is contained
in a block $\{0,y,y'\}$  to which it is connected by an edge
and the block $\{p,q,y\}$ to which it is not connected. There are $N=(n-5)/2$
further blocks containing $y$. Suppose that $h$ of these blocks $\{v_i,w_i,x\}$,
$i= 1,2,\ldots,h$ are not connected to $y$, i.e. $\{y,\{v_i,w_i,y\}\} \notin T$. So without loss
of generality let $\{v_i,\{v_i,w_i,y\}\} \in T$ if $i=1,2,\ldots,h$. Then the $N-h$ blocks
$\{v_i,w_i,y\}$, $i=h+1,h+2,\ldots,N$ are connected to $y$, i.e.
$\{y,\{v_i,w_i,y\}\} \in T$. Clearly $0 \leq h \leq N$ and since $y$ is an odd point,
then $h$ is even. Proceed as follows.\\
(E) If any of the points $v_i$, $i= 1,2,\ldots,h$ are odd points, say $v_\gamma$, then
remove the edge $\{v_\gamma,\{v_\gamma,w_\gamma,y\}\}$
and replace it by $\{y,\{v_\gamma,w_\gamma,y\}\}$. Return to Step 1.\\ 
(F) If any of the points $v_i$ or $w_i$, $i= h+1,h+2,\ldots,N$ are odd points,
again say without loss of generality $v_\gamma$, then
remove the edge $\{y,\{v_\gamma,w_\gamma,y\}\}$ and replace it by $\{v_\gamma,\{v_\gamma,w_\gamma,y\}\}$.
Return to Step 1. 

Step 5. If this step is reached then none of the conditions (A) to (F) are satisfied,
i.e. all of the points $v_i$ and $w_i$ referred to are even points. We count these.
Let $S_p= \{v_i:i=1,2,\ldots,f\} \cup \{v_i,w_i:i=f+1,f+2,\ldots,N\}$\\
and similarly $S_q= \{v_i:i=1,2,\ldots,g\} \cup \{v_i,w_i:i=g+1,g+2,\ldots,N\}$\\         
and $S_y= \{v_i:i=1,2,\ldots,h\} \cup \{v_i,w_i:i=h+1,h+2,\ldots,N\}$.\\
Then $|S_p|=f+2(N-f)=2N-f=(n-5)-f$. Now $0 \leq f \leq N$ and $f$ is even.
Thus if $n \equiv 1$ or 9 (mod 12), $f=(n-5)/2$ and so $|S_p| \geq (n-5)/2$. 
Similarly $|S_q| \geq (n-5)/2$ and $|S_y|\geq (n-5)/2$ so $|S_p|+|S_q|+|S_y| \geq 3(n-5)/2$. 
Now none of the sets $S_p$, $S_q$ nor $S_y$ contain any of the points $0$, $p$, $q$ or $y$.
Since if $n \geq 9$, $3(n-5)/2 > n-4$ at least two of the sets $S_p$, $S_q$ or $S_y$ have
a common point. Further if $n \equiv$ 3 or 7 (mod 12), $f=(n-7)/2$ and so $|S_p| \geq (n-3)/2$. 
Similarly $|S_q| \geq (n-3)/2$ and $|S_y|\geq (n-3)/2$ so $|S_p|+|S_q|+|S_y| \geq 3(n-3)/2$.
Again if $n \geq 3$, $3(n-3)/2 > n-4$ and so at least two of the sets $S_p$, $S_q$ or $S_y$ have
a common point. In either case suppose without loss of generality that
$z \in S_p \cap S_q$.  Denote by $A$ the block which contains the pair $\{p,z\}$
and by $B$ the block which contains the pair $\{q,z\}$. There are four possibilities
to consider.\\
(i) The edges $\{p,A\}$ and $\{q,B\} \in T$.\\Remove both of these edges from $T$
and replace by $\{z,A\}$ and $\{z,B\}$.\\ 
(ii) The edges $\{z,A\}$ and $\{q,B\} \in T$.\\Remove both of these edges from $T$
and replace by $\{p,A\}$ and $\{z,B\}$.\\
(iii) The edges $\{p,A\}$ and $\{z,B\} \in T$.\\Remove both of these edges from $T$
and replace by $\{z,A\}$ and $\{q,B\}$.\\
(iv) The edges $\{z,A\}$ and $\{z,B\} \in T$.\\Remove both of these edges from $T$
and replace by $\{p,A\}$ and $\{q,B\}$.\\
In all four cases return to Step 1. \hfill $\Box$

\section{Latin squares}\label{Latin}
Now we turn our attention to the case of Latin squares of odd order. Let LS($n$) be
a Latin square of odd order. Denote the sets of row points, column points and entry points by
$\mathcal{R}$, $\mathcal{C}$ and $\mathcal{E}$ respectively. Let $\mathcal{B}$ be the set of triples, which
we will also for convenience refer to as blocks, $\{i_r,j_c,k_e\}$ where $i_r \in \mathcal{R}$,
$j_c \in \mathcal{C}$, $k_e \in \mathcal{E}$ and $k=L(i,j)$. The \emph{point-block incidence graph} is
the bipartite graph with vertex set $\mathcal{R} \cup \mathcal{C} \cup \mathcal{E} \cup \mathcal{B}$ and edge set
$\{(v_r,B):v_r \in \mathcal{R}, B \in \mathcal{B}, v_r \in B\} \cup
\{(v_c,B):v_c \in \mathcal{C}, B \in \mathcal{B}, v_c \in B\} \cup
\{(v_e,B):v_e \in \mathcal{E}, B \in \mathcal{B}, v_e \in B\}$. The following result is an exact
analogy of Proposition \ref{spanSTS} above for Steiner triple systems and establishes Theorem \ref{thm:LSodd}. However the proof
is much simpler. We are able to construct the spanning tree with the appropriate property
directly rather than by an iterative procedure.
\begin{proposition}\label{spanLSodd}
Let $L$ = LS($n$) be a Latin square of odd order $n$ and let
$G=G(L)$ be its point-block incidence graph. Then $G$ admits a spanning tree such that every point vertex has even valency in the corresponding co-tree.
\end{proposition}
\textbf{Proof.}
Let $\mathcal{R}=\{0_r,1_r,\ldots,(n-1)_r\}$, $\mathcal{C}=\{0_c,1_c,\ldots,(n-1)_c\}$
and $\mathcal{E}=\{0_e,1_e,\ldots,(n-1)_e\}$. Construct a spanning tree $T$ of the
point-block incidence graph $G$ as follows. For convenience, we will again
refer to point vertices simply as points and block vertices as blocks.  
Let the root (Level 0) of the tree be the point $0_r$. Connect the point to all $n$ blocks
containing it, which will be at Level 1. At Level 2 put all the $2n$ points $x \in \mathcal{C} \cup \mathcal{E}$
and connect each of these to the unique block at Level 1 which contains it.
At Level 3 put the remaining $n-1$ blocks which contain the point $0_c$ and connect
these to $0_c$. Then at Level 4 put the $n-1$ points $\mathcal{R} \setminus \{0_r\}$ and connect each of
these to the unique block at Level 3 which contains it.
In standard form, i.e. $L(0,j)=j$, the structure of the tree is shown below.

\begin{figure}[h]
\begin{center}
      \scalebox{0.45}
       {\includegraphics{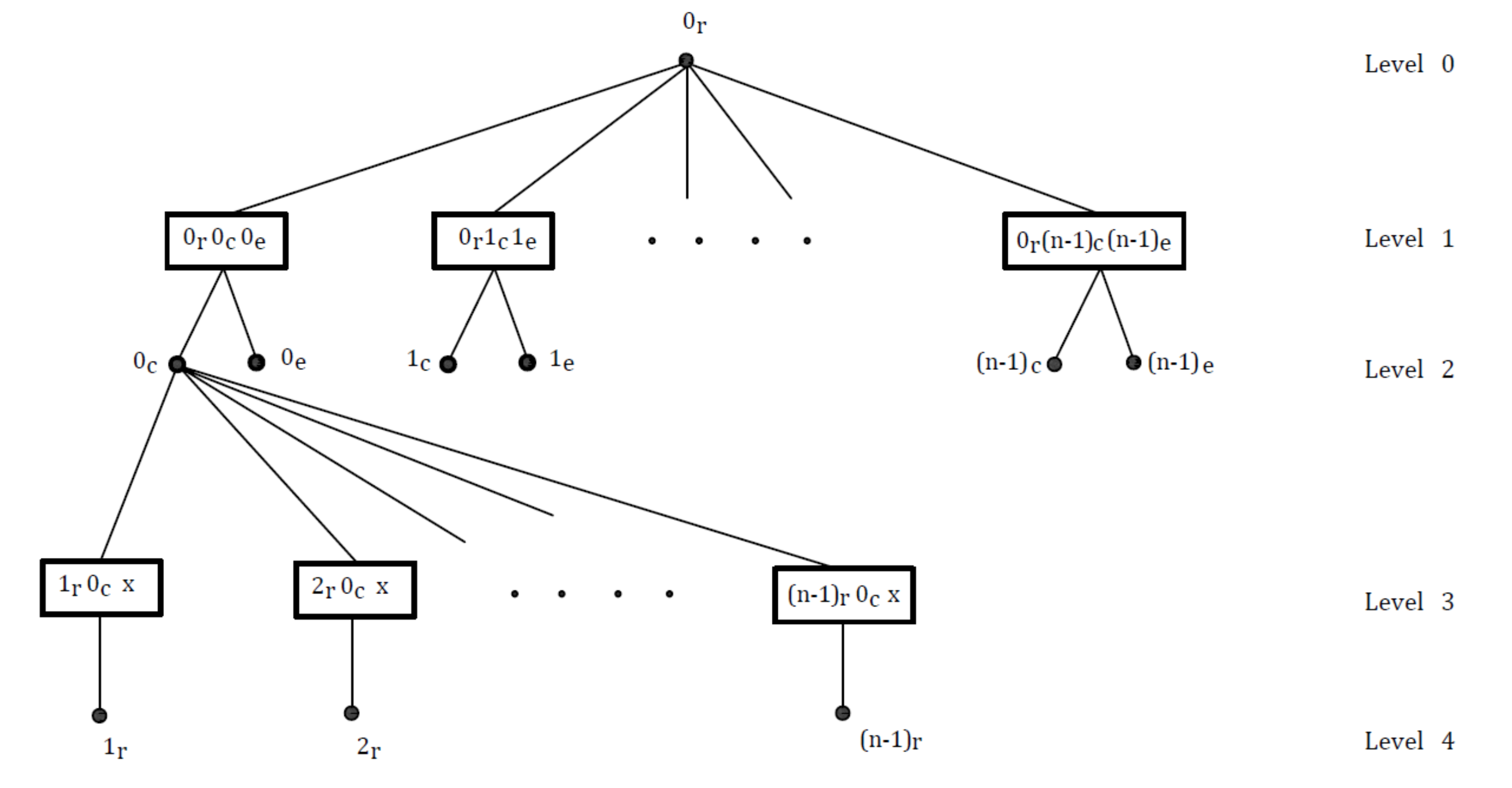}}
       \vskip -2mm
       \caption{Spanning tree of $G(L)$.}
\end{center}
    \end{figure}

 
At this stage we have a tree which contains all the points but only the blocks
which contain the points $0_r$ or $0_c$. Now consider the point $1_c$.
There are $n-1$ blocks containing the point which as yet are not in the tree.
Connect all of these to $1_c$. Repeat this procedure for all of the points in the set
$\mathcal{C} \setminus \{0_c,1_c\}$. We now have a spanning tree of the point-block incidence graph
in which the valency of the points $0_r$ and $x \in \mathcal{C}$ are $n$ and all other points have
valency 1. Thus in the co-tree all points have even valency.
\hfill $\Box$ 

~\\~\\
\noindent 
\Large{\textbf{Acknowledgement}}\\
\normalsize
~\\
The third author acknowledges support from the APVV Research Grants
15-0220 and 17-0428 and the VEGA Research Grants 1/0142/17 and 1/0238/19.

\end{document}